\documentclass[12pt,thmsa]{article}
\usepackage{amssymb}
\usepackage{amsfonts}
\usepackage{sw20bams}

\input tcilatex
\begin{document}

\author{Steven R. Finch}
\title{Rank-$3$ Projections of a $4$-Cube}
\date{May 10, 2012}
\maketitle

\begin{abstract}
The orthogonal projection of a $4$-cube onto a uniform random $3$-subspace
in $\mathbb{R}^{4}$ is a convex $3$-polyhedron $P$ with $14$ vertices almost
surely. Three numerical characteristics of $P$ -- volume, surface area and
mean width -- are studied. These quantities, along with the Euler
characteristic, form a basis of the space of all additive continuous
measures that are invariant under rigid motions in $\mathbb{R}^{3}$. While
computing statistics of $\{vl$, $ar$, $mw\}$, we encounter the generalized
hypergeometric function, elliptic integrals and Catalan's constant. A\ new
constant $7.1185587167...$ also arises and deserves further attention.
\end{abstract}

\footnotetext{%
Copyright \copyright\ 2012 by Steven R. Finch. All rights reserved.}A planar
shadow of a $3$-cube \{$4$-cube\} is a convex hexagon \{octagon\} almost
surely. In an earlier paper \cite{Fi0}, joint moments of hexagonal
\{octagonal\} area and perimeter were computed. It is natural to speculate
on $n$-cubes, under the action of corank-$1$ projections rather than rank-$2$%
.

Let $C$ denote a $4$-cube with edges of unit length, centered at the origin.
To generate a random $3$-subspace $S$ in $\mathbb{R}^{4}$, we select a
random point $U$ uniformly on the $3$-sphere of unit radius. The desired
subspace is the set of all vectors orthogonal to $U$.

We then project the (fixed) $4$-cube $C$ orthogonally onto $S$. This is done
by forming the convex hull of images of all vertices of $C$. The resultant
polyhedron in the hyperplane has $14$ vertices almost surely. More
precisely, if $U=(x,y,z,w)$ is of unit length, then the matrix%
\[
M_{4}=\left( 
\begin{array}{cccc}
\sqrt{1-x^{2}} & -\frac{x\,y}{\sqrt{1-x^{2}}} & -\frac{x\,z}{\sqrt{1-x^{2}}}
& -\frac{x\,w}{\sqrt{1-x^{2}}} \\ 
0 & \sqrt{\frac{z^{2}+w^{2}}{1-x^{2}}} & -\frac{y\,z}{\sqrt{%
(1-x^{2})(z^{2}+w^{2})}} & -\frac{y\,w}{\sqrt{(1-x^{2})(z^{2}+w^{2})}} \\ 
0 & 0 & \frac{w}{\sqrt{z^{2}+w^{2}}} & -\frac{z}{\sqrt{z^{2}+w^{2}}} \\ 
0 & 0 & 0 & 0%
\end{array}%
\right) 
\]%
projects $C$ orthogonally onto a hyperplane, rotated in $\mathbb{R}^{4}$ to
coincide with the $3$-subspace spanned by $(1,0,0,0)$, $(0,1,0,0)$ and $%
(0,0,1,0)$ for convenience. \ Let $T$ be the $4$-row matrix whose columns
constitute all vertices of $C$. Then the first $3$ rows of $M_{4}T$
constitute all images of the vertices in $\mathbb{R}^{3}$ and $3$%
-dimensional convex hull algorithms apply naturally. Such calculations
provide the underpinning for our work.

Spherical coordinates in $\mathbb{R}^{4}$:%
\[
\begin{array}{ccccccc}
x=\cos \theta \sin \varphi \sin \psi , &  & y=\sin \theta \sin \varphi \sin
\psi , &  & z=\cos \varphi \sin \psi , &  & w=\cos \psi%
\end{array}%
\]%
will be used throughout for $U$, where $0\leq \theta <2\pi $, $0\leq \varphi
\leq \pi $, $0\leq \psi \leq \pi $. \ The corresponding Jacobian determinant
is $\sin \varphi \sin ^{2}\psi $; it is best to think of $(\theta ,\varphi
,\psi )$ as possessing joint density $\frac{1}{2\pi ^{2}}\sin \varphi \sin
^{2}\psi $.

\section{Volume}

The projected polyhedral volume $vl$ in $\mathbb{R}^{3}$ is equal to%
\[
\left\vert x\right\vert +\left\vert y\right\vert +\left\vert z\right\vert
+\left\vert w\right\vert , 
\]%
given a unit vector $U=(x,y,z,w)$. It follows that%
\begin{eqnarray*}
\mathbb{E}\left( vl\right) &=&4\cdot 16\dint\limits_{0}^{\pi
/2}\,\dint\limits_{0}^{\pi /2}\,\dint\limits_{0}^{\pi /2}\cos \psi \,\frac{1%
}{2\pi ^{2}}\sin \varphi \sin ^{2}\psi \,d\psi \,d\varphi \,d\theta \\
&=&\frac{16}{3\pi }=1.697652726313550...
\end{eqnarray*}%
\begin{eqnarray*}
\mathbb{E}\left( vl^{2}\right) &=&4\cdot 16\dint\limits_{0}^{\pi
/2}\,\dint\limits_{0}^{\pi /2}\,\dint\limits_{0}^{\pi /2}\cos ^{2}\psi \,%
\frac{1}{2\pi ^{2}}\sin \varphi \sin ^{2}\psi \,d\psi \,d\varphi \,d\theta \\
&&+12\cdot 16\dint\limits_{0}^{\pi /2}\,\dint\limits_{0}^{\pi
/2}\,\dint\limits_{0}^{\pi /2}\cos \varphi \sin \psi \cos \psi \,\frac{1}{%
2\pi ^{2}}\sin \varphi \sin ^{2}\psi \,d\psi \,d\varphi \,d\theta \\
&=&1+\frac{6}{\pi }=2.909859317102744....
\end{eqnarray*}%
More generally, starting with a unit $n$-cube, the projected $(n-1)$%
-polyhedral volume $vl$ in $\mathbb{R}^{n-1}$ satisfies \ 
\[
\mathbb{E}\left( vl\right) =\frac{n}{\sqrt{\pi }}\frac{\Gamma \left( \frac{n%
}{2}\right) }{\Gamma \left( \frac{n+1}{2}\right) }, 
\]%
\[
\mathbb{E}\left( vl^{2}\right) =1+\frac{2(n-1)}{\pi }. 
\]%
These volume moment formulas are the same as formulas in \cite{Fi1} for the
mean width and mean square width, respectively, of the $n$-cube itself. \
Such duality is a special case of a theorem proved in \cite{Du1, Du2, Zn1,
Zn2}. We also have $\max \left( vl\right) =\sqrt{n}$ and $\min \left(
vl\right) =1$.

\section{Surface Area}

The projected polyhedral surface area $ar$ in $\mathbb{R}^{3}$ is equal to%
\[
2\left( \sqrt{x^{2}+y^{2}}+\sqrt{x^{2}+z^{2}}+\sqrt{x^{2}+w^{2}}+\sqrt{%
y^{2}+z^{2}}+\sqrt{y^{2}+w^{2}}+\sqrt{z^{2}+w^{2}}\right) , 
\]%
given a unit vector $U=(x,y,z,w)$. It follows that 
\[
\mathbb{E}\left( ar\right) =6\cdot 16\cdot 2\dint\limits_{0}^{\pi
/2}\,\dint\limits_{0}^{\pi /2}\,\dint\limits_{0}^{\pi /2}\sqrt{\cos
^{2}\varphi \sin ^{2}\psi +\cos ^{2}\psi }\,\frac{1}{2\pi ^{2}}\sin \varphi
\sin ^{2}\psi \,d\psi \,d\varphi \,d\theta =8 
\]%
and, since $x^{2}+y^{2}=\sin ^{2}\varphi \sin ^{2}\psi $, 
\begin{eqnarray*}
\mathbb{E}\left( ar^{2}\right) &=&6\cdot 16\cdot 4\dint\limits_{0}^{\pi
/2}\,\dint\limits_{0}^{\pi /2}\,\dint\limits_{0}^{\pi /2}\left( \cos
^{2}\varphi \sin ^{2}\psi +\cos ^{2}\psi \right) \,\frac{1}{2\pi ^{2}}\sin
\varphi \sin ^{2}\psi \,d\psi \,d\varphi \,d\theta \\
&&+24\cdot 16\cdot 4\dint\limits_{0}^{\pi /2}\,\dint\limits_{0}^{\pi
/2}\,\dint\limits_{0}^{\pi /2}\sin \varphi \sin \psi \sqrt{\sin ^{2}\theta
\sin ^{2}\varphi \sin ^{2}\psi +\cos ^{2}\psi }\,\frac{1}{2\pi ^{2}}\sin
\varphi \sin ^{2}\psi \,d\psi \,d\varphi \,d\theta \\
&&+6\cdot 16\cdot 4\dint\limits_{0}^{\pi /2}\,\dint\limits_{0}^{\pi
/2}\,\dint\limits_{0}^{\pi /2}\sin \varphi \sin \psi \sqrt{\cos ^{2}\varphi
\sin ^{2}\psi +\cos ^{2}\psi }\,\frac{1}{2\pi ^{2}}\sin \varphi \sin
^{2}\psi \,d\psi \,d\varphi \,d\theta \\
&=&12+6\zeta _{4}+3\pi =64.136130261087789....
\end{eqnarray*}%
Details on derivation of the mean square result, including a formula for the
constant $\zeta _{4}$, will appear shortly. We also have $\max \left(
ar\right) =6\sqrt{2}$ and $\min \left( ar\right) =6$.

More generally, starting with a unit $n$-cube, the projected $(n-1)$%
-polyhedral surface area $ar$ in $\mathbb{R}^{n-1}$ satisfies%
\[
\mathbb{E}\left( ar\right) =\frac{\sqrt{\pi }(n-1)n}{2}\frac{\Gamma \left( 
\frac{n}{2}\right) }{\Gamma \left( \frac{n+1}{2}\right) }, 
\]%
\[
\mathbb{E}\left( ar^{2}\right) =4(n-1)+(n-2)(n-1)\zeta _{n}+\frac{%
(n-3)(n-2)(n-1)}{2}\pi . 
\]%
For the special case $n=3$, we have%
\[
\frac{ar}{2}=\sqrt{x^{2}+y^{2}}+\sqrt{x^{2}+z^{2}}+\sqrt{y^{2}+z^{2}}=\sqrt{%
1-x^{2}}+\sqrt{1-y^{2}}+\sqrt{1-z^{2}}=\frac{\pi }{2}\,mw 
\]%
given that $x^{2}+y^{2}+z^{2}=1$, which implies a closed-form expression for 
$\zeta _{3}$ (using results from an upcoming subsection \textquotedblleft
2D\ Analog\textquotedblright ). Numerical evidence strongly suggests that $%
\zeta _{n}=\zeta _{3}$ for all $n\geq 4$, but a rigorous proof is not known.

\section{Mean Width}

The projected polyhedral mean width $mw$ in $\mathbb{R}^{3}$ is equal to 
\[
\frac{1}{2}\left( \sqrt{1-x^{2}}+\sqrt{1-y^{2}}+\sqrt{1-z^{2}}+\sqrt{1-w^{2}}%
\right) , 
\]%
given a unit vector $U=(x,y,z,w)$. It follows that%
\begin{eqnarray*}
\mathbb{E}\left( mw\right) &=&4\cdot 16\cdot \frac{1}{2}\dint\limits_{0}^{%
\pi /2}\,\dint\limits_{0}^{\pi /2}\,\dint\limits_{0}^{\pi /2}\sqrt{1-\cos
^{2}\psi }\,\frac{1}{2\pi ^{2}}\sin \varphi \sin ^{2}\psi \,d\psi \,d\varphi
\,d\theta \\
&=&\frac{16}{3\pi }=\mathbb{E}\left( vl\right) ,
\end{eqnarray*}%
\begin{eqnarray*}
\mathbb{E}\left( mw^{2}\right) &=&4\cdot 16\cdot \frac{1}{4}%
\dint\limits_{0}^{\pi /2}\,\dint\limits_{0}^{\pi /2}\,\dint\limits_{0}^{\pi
/2}\left( 1-\cos ^{2}\psi \right) \frac{1}{2\pi ^{2}}\sin \varphi \sin
^{2}\psi \,d\psi \,d\varphi \,d\theta \\
&&+12\cdot 16\cdot \frac{1}{4}\dint\limits_{0}^{\pi
/2}\,\dint\limits_{0}^{\pi /2}\,\dint\limits_{0}^{\pi /2}\sqrt{1-\cos
^{2}\varphi \sin ^{2}\psi }\sqrt{1-\cos ^{2}\psi }\,\frac{1}{2\pi ^{2}}\sin
\varphi \sin ^{2}\psi \,d\psi \,d\varphi \,d\theta \\
&=&3\left( \frac{1}{4}+\frac{\pi }{8}+\frac{1}{\pi }\right)
=2.883026903647544...<1+\frac{6}{\pi }=\mathbb{E}\left( vl^{2}\right) \text{.%
}
\end{eqnarray*}%
We also have $\max \left( mw\right) =\sqrt{3}$ and $\min \left( mw\right)
=3/2$.

More generally, starting with a unit $n$-cube, the projected $(n-1)$%
-polyhedral mean width $mw$ in $\mathbb{R}^{n-1}$ satisfies \ 
\[
\mathbb{E}\left( mw\right) =\frac{n}{\sqrt{\pi }}\frac{\Gamma \left( \frac{n%
}{2}\right) }{\Gamma \left( \frac{n+1}{2}\right) }=\mathbb{E}\left(
vl\right) 
\]%
but a general formula for $\mathbb{E}\left( mw^{2}\right) $ is unknown. We
examine the cases $n=3$ and $n=5$ in the following subsections.\pagebreak

\subsection{2D\ Analog}

As outlined in \cite{Fi0},

\[
mw=\frac{2}{\pi }\left( \sqrt{1-x^{2}}+\sqrt{1-y^{2}}+\sqrt{1-z^{2}}\right) 
\]%
given that $x^{2}+y^{2}+z^{2}=1$, hence $\mathbb{E}\left( mw\right) =3/2$ and

\begin{eqnarray*}
\mathbb{E}\left( mw^{2}\right) &=&3\cdot 8\cdot \left( \frac{2}{\pi }\right)
^{2}\,\dint\limits_{0}^{\pi /2}\,\dint\limits_{0}^{\pi /2}\left( 1-\cos
^{2}\varphi \right) \frac{1}{4\pi }\sin \varphi \,d\varphi \,d\theta \\
&&+6\cdot 8\cdot \left( \frac{2}{\pi }\right) ^{2}\dint\limits_{0}^{\pi
/2}\,\dint\limits_{0}^{\pi /2}\,\sqrt{1-\cos ^{2}\theta \sin ^{2}\varphi }%
\sqrt{1-\cos ^{2}\varphi }\,\frac{1}{4\pi }\sin \varphi \,d\varphi \,d\theta
\\
&=&\frac{2}{\pi ^{2}}\left[ 4+3\pi \,_{3}F_{2}\left( -\tfrac{1}{2},\tfrac{1}{%
2},\tfrac{3}{2};1,2;1\right) \right] =2.253091059149751...<1\mathbb{+}\frac{4%
}{\pi }
\end{eqnarray*}%
where 
\[
_{3}F_{2}(a_{1},a_{2},a_{3};b_{1},b_{2};z)=\frac{\Gamma (b_{1})\Gamma (b_{2})%
}{\Gamma (a_{1})\Gamma (a_{2})\Gamma (a_{3})}\dsum\limits_{k=0}^{\infty }%
\frac{\Gamma (a_{1}+k)\Gamma (a_{2}+k)\Gamma (a_{3}+k)}{\Gamma
(b_{1}+k)\Gamma (b_{2}+k)}\frac{z^{k}}{k!} 
\]%
is the generalized hypergeometric function. As a byproduct, the constant $%
\zeta _{3}$ defined earlier is equal to $3\pi \,_{3}F_{2}\left( -\tfrac{1}{2}%
,\tfrac{1}{2},\tfrac{3}{2};1,2;1\right) $. We also have $\max \left(
mw\right) =2\sqrt{6}/\pi $ and $\min \left( mw\right) =4/\pi $.\pagebreak\ \ 

\subsection{4D\ Analog}

We here have%
\[
mw=\frac{4}{3\pi }\dsum\limits_{j=1}^{5}\sqrt{1-x_{j}^{2}} 
\]%
given that $\sum_{j=1}^{5}x_{j}^{2}=1$, hence $\mathbb{E}\left( mw\right)
=15/8$ and%
\[
\mathbb{E}\left( mw^{2}\right) =32\cdot \left( \frac{4}{3\pi }\right)
^{2}\left( 5I+20J\right) 
\]%
where%
\[
I=\dint\limits_{0}^{\pi /2}\,\dint\limits_{0}^{\pi
/2}\,\dint\limits_{0}^{\pi /2}\,\dint\limits_{0}^{\pi /2}\left( 1-\cos
^{2}\varphi _{3}\right) \frac{3}{8\pi ^{2}}\sin \varphi _{1}\sin ^{2}\varphi
_{2}\sin ^{3}\varphi _{3}\,d\varphi _{3}d\varphi _{2}d\varphi _{1}d\theta , 
\]%
\[
J=\dint\limits_{0}^{\pi /2\,}\dint\limits_{0}^{\pi
/2}\,\dint\limits_{0}^{\pi /2}\,\dint\limits_{0}^{\pi /2}\sqrt{1-\cos
^{2}\varphi _{2}\sin ^{2}\varphi _{3}}\sqrt{1-\cos ^{2}\varphi _{3}}\,\frac{3%
}{8\pi ^{2}}\sin \varphi _{1}\sin ^{2}\varphi _{2}\sin ^{3}\varphi
_{3}\,d\varphi _{3}d\varphi _{2}d\varphi _{1}d\theta . 
\]%
The expression for $\mathbb{E}\left( mw^{2}\right) $ simplifies to

\begin{eqnarray*}
&&\frac{4}{81\pi ^{4}}\left[ 144\pi ^{2}-10\Gamma \left( \tfrac{1}{4}\right)
^{4}+45\pi ^{3}\left( 8\,_{3}F_{2}\left( -\tfrac{1}{2},\tfrac{1}{2},\tfrac{3%
}{2};1,2;1\right) -\,_{3}F_{2}\left( -\tfrac{1}{2},\tfrac{1}{2},\tfrac{3}{2}%
;1,3;1\right) \right) \right] \\
&=&3.516040901689803...<1\mathbb{+}\frac{8}{\pi }.
\end{eqnarray*}%
We also have $\max \left( mw\right) =8\sqrt{5}/(3\pi )$ and $\min \left(
mw\right) =16/(3\pi )$. \ The fact that $\mathbb{E}\left( mw\right) $ at one
level becomes $\min \left( mw\right) $ at the next level is
interesting.\pagebreak

\section{Details}

With regard to $\mathbb{E}\left( ar^{2}\right) $, we first prove that%
\[
\dint\limits_{0}^{\pi /2}\,\dint\limits_{0}^{\pi /2}\,\dint\limits_{0}^{\pi
/2}\sin \varphi \sin \psi \sqrt{\cos ^{2}\varphi \sin ^{2}\psi +\cos
^{2}\psi }\,\frac{1}{2\pi ^{2}}\sin \varphi \sin ^{2}\psi \,d\psi \,d\varphi
\,d\theta =\frac{\pi }{128}. 
\]%
After rearranging and integrating out $\theta $, the integral becomes 
\begin{eqnarray*}
&&\frac{\pi }{2}\dint\limits_{0}^{\pi /2}\,\cos \psi \sin ^{3}\psi
\dint\limits_{0}^{\pi /2}\sqrt{\cos ^{2}\varphi \tan ^{2}\psi +1}\,\frac{1}{%
2\pi ^{2}}\sin ^{2}\varphi \,d\varphi \,d\psi \\
&=&\frac{1}{4\pi }\dint\limits_{0}^{\pi /2}\,\cos \psi \sin ^{3}\psi \frac{%
\left( -1+\tan ^{2}\psi \right) E\left( i\,\tan \psi \right) +\left( 1+\tan
^{2}\psi \right) K\left( i\,\tan \psi \right) }{3\tan ^{2}\psi }d\psi
\end{eqnarray*}%
where $i$ is the imaginary unit and%
\[
\begin{array}{ccc}
K(\xi )=\dint\limits_{0}^{\pi /2}\dfrac{1}{\sqrt{1-\xi ^{2}\sin (\theta )^{2}%
}}\,d\theta , &  & E(\xi )=\dint\limits_{0}^{\pi /2}\sqrt{1-\xi ^{2}\sin
(\theta )^{2}}\,d\theta%
\end{array}%
\]%
are complete elliptic integrals of the first and second kind. Let $t=\tan
\psi $, then%
\[
dt=\sec ^{2}\psi \,d\psi =\left( 1+t^{2}\right) d\psi 
\]%
hence%
\begin{eqnarray*}
\frac{\cos \psi \sin ^{3}\psi }{\tan ^{2}\psi }d\psi &=&\cos ^{3}\psi \sin
\psi \,d\psi =\left( \frac{1}{1+t^{2}}\right) ^{3/2}\left( 1-\frac{1}{1+t^{2}%
}\right) ^{1/2}d\psi \\
&=&\left( \frac{1}{1+t^{2}}\right) ^{5/2}\left( 1-\frac{1}{1+t^{2}}\right)
^{1/2}dt \\
&=&\frac{t}{\left( 1+t^{2}\right) ^{3}}dt
\end{eqnarray*}%
and the integral becomes%
\begin{eqnarray*}
&&\frac{1}{12\pi }\dint\limits_{0}^{\infty }\,\frac{\left( -1+t^{2}\right)
E\left( i\,t\right) +\left( 1+t^{2}\right) K\left( i\,t\right) }{\left(
1+t^{2}\right) ^{3}}t\,dt \\
&=&\frac{1}{12\pi }\left[ \dint\limits_{0}^{\infty }\,\frac{E\left(
i\,t\right) }{\left( 1+t^{2}\right) ^{2}}t\,dt-2\dint\limits_{0}^{\infty }\,%
\frac{E\left( i\,t\right) }{\left( 1+t^{2}\right) ^{3}}t\,dt+\dint%
\limits_{0}^{\infty }\,\frac{K\left( i\,t\right) }{\left( 1+t^{2}\right) ^{2}%
}t\,dt\right] \\
&=&\frac{\pi }{96}-\frac{\pi }{128}+\frac{\pi }{192}=\frac{\pi }{128}.
\end{eqnarray*}

The other integral%
\[
\dint\limits_{0}^{\pi /2}\,\dint\limits_{0}^{\pi /2}\,\dint\limits_{0}^{\pi
/2}\sin \varphi \sin \psi \sqrt{\sin ^{2}\theta \sin ^{2}\varphi \sin
^{2}\psi +\cos ^{2}\psi }\,\frac{1}{2\pi ^{2}}\sin \varphi \sin ^{2}\psi
\,d\psi \,d\varphi \,d\theta 
\]%
is harder. We rearrange it as%
\begin{eqnarray*}
&&\dint\limits_{0}^{\pi /2}\,\dint\limits_{0}^{\pi /2}\cos \psi \sin
^{3}\psi \dint\limits_{0}^{\pi /2}\sqrt{\sin ^{2}\theta \sin ^{2}\varphi
\tan ^{2}\psi +1}\,\frac{1}{2\pi ^{2}}\sin ^{2}\varphi \,d\theta \,d\varphi
\,d\psi  \\
&=&\frac{1}{2\pi ^{2}}\dint\limits_{0}^{\pi /2}\,\cos \psi \sin ^{3}\psi
\dint\limits_{0}^{\pi /2}E\left( i\sin \varphi \tan \psi \right) \sin
^{2}\varphi \,d\varphi \,d\psi  \\
&=&\frac{1}{2\pi ^{2}}\dint\limits_{0}^{\pi /2}\,\cos \psi \sin ^{3}\psi \,%
\frac{f(\psi )+g(\psi )+h(\psi )}{3\tan ^{2}\psi }\,d\psi ,
\end{eqnarray*}%
where 
\[
f(\psi )=\left( -2+4\tan ^{2}\psi \right) E\left( i\,\sqrt{\tfrac{\sec \psi
-1}{2}}\right) ^{2},
\]%
\[
g(\psi )=2\left( 1-2\tan ^{2}\psi +\sec \psi \right) E\left( i\,\sqrt{\tfrac{%
\sec \psi -1}{2}}\right) K\left( i\,\sqrt{\tfrac{\sec \psi -1}{2}}\right) ,
\]%
\[
h(\psi )=\left( -1+\tan ^{2}\psi -\sec \psi \right) K\left( i\,\sqrt{\tfrac{%
\sec \psi -1}{2}}\right) ^{2}.
\]%
Let $t=\sqrt{(\sec \psi -1)/2}$, then $1+2t^{2}=\sec \psi $ and%
\[
4t\,dt=\sec \psi \tan \psi \,d\psi =\left( 1+2t^{2}\right) \left( \left(
1+2t^{2}\right) ^{2}-1\right) ^{1/2}d\psi 
\]%
hence 
\begin{eqnarray*}
\frac{\cos \psi \sin ^{3}\psi }{\tan ^{2}\psi }d\psi  &=&\cos ^{3}\psi \sin
\psi \,d\psi =\left( \frac{1}{1+2t^{2}}\right) ^{3}\left( 1-\frac{1}{\left(
1+2t^{2}\right) ^{2}}\right) ^{1/2}d\psi  \\
&=&(4t)\left( \frac{1}{1+2t^{2}}\right) ^{4}\left( 1-\frac{1}{\left(
1+2t^{2}\right) ^{2}}\right) ^{1/2}\frac{1}{\left( \left( 1+2t^{2}\right)
^{2}-1\right) ^{1/2}}dt \\
&=&\frac{4t}{\left( 1+2t^{2}\right) ^{5}}dt.
\end{eqnarray*}%
The integral becomes%
\begin{eqnarray*}
&&\frac{2}{3\pi ^{2}}\dint\limits_{0}^{\infty }\,\frac{\left[ -2+4\left(
\left( 1+2t^{2}\right) ^{2}-1\right) \right] E\left( i\,t\right) ^{2}}{%
\left( 1+2t^{2}\right) ^{5}}t\,dt \\
&&+\frac{4}{3\pi ^{2}}\dint\limits_{0}^{\infty }\,\frac{\left[ 1-2\left(
\left( 1+2t^{2}\right) ^{2}-1\right) +(1+2t^{2})\right] E\left( i\,t\right)
K\left( i\,t\right) }{\left( 1+2t^{2}\right) ^{5}}t\,dt \\
&&+\frac{2}{3\pi ^{2}}\dint\limits_{0}^{\infty }\,\frac{\left[ -1+\left(
\left( 1+2t^{2}\right) ^{2}-1\right) -(1+2t^{2})\right] K\left( i\,t\right)
^{2}}{\left( 1+2t^{2}\right) ^{5}}t\,dt
\end{eqnarray*}%
which simplifies to%
\[
\frac{4}{3\pi ^{2}}\dint\limits_{0}^{\infty }\,\frac{\left(
8t^{4}+8t^{2}-1\right) E\left( i\,t\right) ^{2}-2\left(
4t^{4}+3t^{2}-1\right) E\left( i\,t\right) K\left( i\,t\right) +\left(
2t^{4}+t^{2}-1\right) K\left( i\,t\right) ^{2}}{\left( 2t^{2}+1\right) ^{5}}%
t\,dt.
\]%
Call this expression $\zeta _{4}/256$. The constant $\zeta
_{4}=7.118558716719735...$ coincides with $\zeta _{3}$ to high numerical
precision, but algebraic confirmation remains open.

\section{Correlations}

The joint moment for volume and surface area is 
\begin{eqnarray*}
\mathbb{E}\left( vl\cdot ar\right) &=&12\cdot 16\cdot 2\dint\limits_{0}^{\pi
/2}\,\dint\limits_{0}^{\pi /2}\,\dint\limits_{0}^{\pi /2}\left[ \cos \psi 
\sqrt{\cos ^{2}\varphi \sin ^{2}\psi +\cos ^{2}\psi }\,\right. \\
&&\left. +\cos \psi \sqrt{\sin ^{2}\theta \sin ^{2}\varphi \sin ^{2}\psi
+\cos ^{2}\varphi \sin ^{2}\psi }\,\right] \frac{1}{2\pi ^{2}}\sin \varphi
\sin ^{2}\psi \,d\psi \,d\varphi \,d\theta \\
&=&6\left( 1+\frac{4}{\pi }\right) =13.639437268410976...
\end{eqnarray*}%
which, together with earlier results, implies that the correlation between $%
vl$ and $ar$ is $0.945...$. Likewise, we have%
\begin{eqnarray*}
\mathbb{E}\left( vl\cdot mw\right) &=&4\cdot 16\cdot \frac{1}{2}%
\dint\limits_{0}^{\pi /2}\,\dint\limits_{0}^{\pi /2}\,\dint\limits_{0}^{\pi
/2}\cos \psi \sqrt{1-\cos ^{2}\psi }\,\frac{1}{2\pi ^{2}}\sin \varphi \sin
^{2}\psi \,d\psi \,d\varphi \,d\theta \\
&&+12\cdot 16\cdot \frac{1}{2}\dint\limits_{0}^{\pi
/2}\,\dint\limits_{0}^{\pi /2}\,\dint\limits_{0}^{\pi /2}\cos \varphi \sin
\psi \sqrt{1-\cos ^{2}\psi }\,\frac{1}{2\pi ^{2}}\sin \varphi \sin ^{2}\psi
\,d\psi \,d\varphi \,d\theta \\
&=&\frac{9}{4}+\frac{2}{\pi }=2.886619772367581...,
\end{eqnarray*}%
\begin{eqnarray*}
\mathbb{E}\left( ar\cdot mw\right) &=&12\cdot 16\cdot 2\cdot \frac{1}{2}%
\dint\limits_{0}^{\pi /2}\,\dint\limits_{0}^{\pi /2}\,\dint\limits_{0}^{\pi
/2}\,\left[ \sqrt{\cos ^{2}\varphi \sin ^{2}\psi +\cos ^{2}\psi }\sqrt{%
1-\cos ^{2}\psi }\right. \\
&&\left. +\sqrt{\sin ^{2}\theta \sin ^{2}\varphi \sin ^{2}\psi +\cos
^{2}\varphi \sin ^{2}\psi }\sqrt{1-\cos ^{2}\psi }\,\right] \frac{1}{2\pi
^{2}}\sin \varphi \sin ^{2}\psi \,d\psi \,d\varphi \,d\theta \\
&=&\frac{3\left( 5+2G\right) }{\pi }+\frac{9\pi }{4}=13.592597187518807...
\end{eqnarray*}%
where 
\[
G=\dsum\limits_{k=0}^{\infty }\frac{(-1)^{k}}{(2k+1)^{2}} 
\]%
is Catalan's constant \cite{Fc}. The correlation between $vl$ and $mw$ is
consequently $0.870...$ and the correlation between $ar$ and $mw$ is $%
0.973...$.

\section{Related Work}

Let us return to an open issue in \cite{Fi0} surrounding rank-$2$
projections of a $4$-cube. On the one hand, octagonal perimeter is equal to%
\[
2\left( \sqrt{1-p^{2}-x^{2}}+\sqrt{1-q^{2}-y^{2}}+\sqrt{1-r^{2}-z^{2}}+\sqrt{%
1-s^{2}-w^{2}}\right) , 
\]%
given orthogonal unit vectors $U=(x,y,z,w)$ and $V=(p,q,r,s)$. \ This
formula makes accurate calculation of the second moment of perimeter $%
28.495...$ more feasible. \ The maximum value of perimeter is $4\sqrt{2}$;
the minimum value is $4$.

On the other hand, octagonal area does not appear to have so simple a
description. Let 
\[
\left( 
\begin{array}{c}
p \\ 
q \\ 
r \\ 
s%
\end{array}%
\right) =\cos \kappa \sin \lambda \left( 
\begin{array}{c}
-y \\ 
x \\ 
-w \\ 
z%
\end{array}%
\right) +\sin \kappa \sin \lambda \left( 
\begin{array}{c}
-z \\ 
w \\ 
x \\ 
-y%
\end{array}%
\right) +\cos \lambda \left( 
\begin{array}{c}
-w \\ 
-z \\ 
y \\ 
x%
\end{array}%
\right) 
\]%
where $0\leq \kappa <2\pi $, $0\leq \lambda \leq \pi $. Let%
\[
a_{1}=r\,y+s\,y-q\,z-s\,z-q\,w+r\,w, 
\]%
\[
a_{2}=r\,y-s\,y-q\,z-s\,z+q\,w+r\,w, 
\]%
\[
a_{3}=r\,y+s\,y-q\,z+s\,z-q\,w-r\,w, 
\]%
\[
b_{1}=p\,q-p\,s+x\,y-x\,w, 
\]%
\[
b_{2}=p\,q-p\,r+x\,y-x\,z, 
\]%
\[
b_{3}=p\,q+p\,s+x\,y+x\,w, 
\]%
\[
c=2(1-p^{2}-x^{2}). 
\]%
Then, in a neighborhood of $(\theta _{0},\varphi _{0},\psi _{0},\kappa
_{0},\lambda _{0})=(1,1,1,1,1)$, octagonal area is equal to%
\[
\frac{1}{c}\left[ a_{1}b_{2}+a_{2}(c-b_{1}+b_{2})+a_{3}b_{1}\right] ; 
\]%
in a neighborhood of $(\theta _{0},\varphi _{0},\psi _{0},\kappa
_{0},\lambda _{0})=\left( \frac{1}{2},\frac{1}{2},\frac{1}{2},\frac{1}{2},%
\frac{1}{2}\right) $, octagonal area is equal to%
\[
\frac{1}{c}\left[ a_{1}b_{2}-a_{2}(b_{1}-b_{2})+a_{3}(c+b_{1})\right] ; 
\]%
in a neighborhood of $(\theta _{0},\varphi _{0},\psi _{0},\kappa
_{0},\lambda _{0})=\left( \frac{3}{4},\frac{3}{4},\frac{3}{4},\frac{3}{4},%
\frac{3}{4}\right) $, octagonal area is equal to%
\[
\frac{1}{c}\left[ a_{1}b_{2}-a_{2}(c+b_{1}-b_{2})+a_{3}b_{1}\right] ; 
\]%
in a neighborhood of $(\theta _{0},\varphi _{0},\psi _{0},\kappa
_{0},\lambda _{0})=\left( \frac{5}{6},\frac{5}{6},\frac{5}{6},\frac{5}{6},%
\frac{5}{6}\right) $, octagonal area is equal to%
\[
\frac{1}{c}\left[ -a_{1}(c-b_{2})-a_{2}(b_{1}-b_{2})+a_{3}b_{1}\right] ; 
\]%
in a neighborhood of $(\theta _{0},\varphi _{0},\psi _{0},\kappa
_{0},\lambda _{0})=\left( \frac{7}{8},\frac{7}{8},\frac{7}{8},\frac{7}{8},%
\frac{7}{8}\right) $, octagonal area is equal to

\[
\frac{1}{c}\left[ a_{1}b_{2}-a_{2}(b_{1}-b_{2})-a_{3}(c-b_{1})\right] ; 
\]%
in a neighborhood of $(\theta _{0},\varphi _{0},\psi _{0},\kappa
_{0},\lambda _{0})=\left( \frac{4}{5},1,\frac{2}{5},\frac{3}{5},\frac{1}{5}%
\right) $, octagonal area is equal to

\[
\frac{1}{c}\left[ a_{1}(c+b_{2})-a_{2}(b_{1}-b_{2})+a_{3}b_{1}\right] ; 
\]%
Other branches also exist, but for reasons of space, we stop here. The
maximum value of octagonal area is $1+\sqrt{2}$ \cite{Zn1, Zn2}; the minimum
value is $1$.

An addendum to \cite{VB} now exists, with voluminous detail on the
derivation of a relevant density function \cite{Vic}.

The present work is the fifth (and final) in a pentalogy that began with 
\cite{Fi2}, continued with \cite{Fi1, Fi3} and then again with \cite{Fi0}.

\section{Acknowledgements}

Wouter Meeussen's package ConvexHull3D.m was helpful to me in preparing this
paper \cite{Msn}. \ He kindly extended the software functionality at my
request. Much more relevant material can be found at \cite{Fi4}, including
experimental computer runs that aided theoretical discussion here.

\section{Addendum}

Here is a proof that $\zeta _{5}=\zeta _{4}$. For $n=5$, we have%
\[
ar=2\dsum\limits_{1\leq j<k\leq 5}\sqrt{x_{j}^{2}+x_{k}^{2}}
\]%
given that $\sum_{j=1}^{5}x_{j}^{2}=1$. The hard integral for $n=4$ in a
preceding subsection \textquotedblleft Details\textquotedblright\ here
becomes 
\begin{eqnarray*}
&&\dint\limits_{0}^{\pi /2}\,\dint\limits_{0}^{\pi
/2}\,\dint\limits_{0}^{\pi /2}\,\dint\limits_{0}^{\pi /2}\sin \varphi
_{1}\sin \varphi _{2}\sin \varphi _{3}\sqrt{\sin ^{2}\theta \sin ^{2}\varphi
_{1}\sin ^{2}\varphi _{2}\sin ^{2}\varphi _{3}+\cos ^{2}\varphi _{3}} \\
&&\;\;\;\;\;\;\;\;\;\;\;\;\;\;\;\;\times \frac{3}{8\pi ^{2}}\sin \varphi
_{1}\sin ^{2}\varphi _{2}\sin ^{3}\varphi _{3}\,d\varphi _{3}d\varphi
_{2}d\varphi _{1}d\theta  \\
&=&\dint\limits_{0}^{\pi /2}\,\dint\limits_{0}^{\pi
/2}\,\dint\limits_{0}^{\pi /2}\,\cos \varphi _{3}\sin ^{4}\varphi
_{3}\dint\limits_{0}^{\pi /2}\sqrt{\sin ^{2}\theta \sin ^{2}\varphi _{1}\sin
^{2}\varphi _{2}\tan ^{2}\varphi _{3}+1} \\
&&\;\;\;\;\;\;\;\;\;\;\;\;\;\;\;\;\times \frac{3}{8\pi ^{2}}\sin ^{2}\varphi
_{1}\sin ^{3}\varphi _{2}\,d\theta \,d\varphi _{1}d\varphi _{2}d\varphi _{3}
\\
&=&\frac{3}{8\pi ^{2}}\dint\limits_{0}^{\pi /2}\,\dint\limits_{0}^{\pi
/2}\,\,\cos \varphi _{3}\sin ^{4}\varphi _{3}\dint\limits_{0}^{\pi
/2}E\left( i\sin \varphi _{1}\sin \varphi _{2}\tan \varphi _{3}\right) \sin
^{2}\varphi _{1}\sin ^{3}\varphi _{2}\,d\varphi _{1}d\varphi _{2}d\varphi
_{3} \\
&=&\frac{3}{8\pi ^{2}}\dint\limits_{0}^{\pi /2}\,\dint\limits_{0}^{\pi
/2}\,\,\sin \varphi _{2}\cos \varphi _{3}\sin ^{4}\varphi _{3}\frac{%
f(\varphi _{2},\varphi _{3})+g(\varphi _{2},\varphi _{3})+h(\varphi
_{2},\varphi _{3})}{3\tan ^{2}\varphi _{3}}d\varphi _{2}d\varphi _{3}
\end{eqnarray*}%
where%
\[
f(\varphi ,\psi )=\left( -2+4\sin ^{2}\varphi \tan ^{2}\psi \right) E\left(
i\,\sqrt{\tfrac{\sqrt{1+\sin ^{2}\varphi \tan ^{2}\psi }-1}{2}}\right) ^{2},
\]%
\begin{eqnarray*}
g(\varphi ,\psi ) &=&2\left( 1-2\sin ^{2}\varphi \tan ^{2}\psi +\sqrt{1+\sin
^{2}\varphi \tan ^{2}\psi }\right)  \\
&&\;\;\;\;\times E\left( i\,\sqrt{\tfrac{\sqrt{1+\sin ^{2}\varphi \tan
^{2}\psi }-1}{2}}\right) K\left( i\,\sqrt{\tfrac{\sqrt{1+\sin ^{2}\varphi
\tan ^{2}\psi }-1}{2}}\right) ,
\end{eqnarray*}%
\[
h(\varphi ,\psi )=\left( -1+\sin ^{2}\varphi \tan ^{2}\psi -\sqrt{1+\sin
^{2}\varphi \tan ^{2}\psi }\right) K\left( i\,\sqrt{\tfrac{\sqrt{1+\sin
^{2}\varphi \tan ^{2}\psi }-1}{2}}\right) ^{2}.
\]%
Let%
\[
\begin{array}{ccc}
u=\sqrt{\tfrac{\sqrt{1+\sin ^{2}\varphi \tan ^{2}\psi }-1}{2}}, &  & v=\sin
\varphi 
\end{array}%
\]%
then%
\[
\begin{array}{ccc}
1+2u^{2}=\sqrt{1+\sin ^{2}\varphi \tan ^{2}\psi }, &  & \sqrt{1-v^{2}}=\cos
\varphi 
\end{array}%
\]%
hence%
\[
\tan ^{2}\psi =\dfrac{\left( 1+2u^{2}\right) ^{2}-1}{v^{2}}=\dfrac{%
4u^{2}\left( 1+u^{2}\right) }{v^{2}},
\]%
\[
\sec ^{2}\psi =1+\dfrac{4u^{2}\left( 1+u^{2}\right) }{v^{2}}=\frac{%
4u^{2}\left( 1+u^{2}\right) +v^{2}}{v^{2}}
\]%
and%
\begin{eqnarray*}
2\left( 1+2u^{2}\right) \left( 4u\right) \frac{\partial u}{\partial \psi }
&=&\left( \sin ^{2}\varphi \right) \left( 2\tan \psi \sec ^{2}\psi \right) 
\\
&=&v^{2}\cdot 2\sqrt{\dfrac{4u^{2}\left( 1+u^{2}\right) }{v^{2}}}\,\frac{%
4u^{2}\left( 1+u^{2}\right) +v^{2}}{v^{2}} \\
&=&\frac{\left( 4u\right) \left[ 4u^{2}\left( 1+u^{2}\right) +v^{2}\right] 
\sqrt{1+u^{2}}}{v}.
\end{eqnarray*}%
It follows that%
\[
\begin{array}{ccc}
\dfrac{\partial u}{\partial \psi }=\dfrac{\left[ 4u^{2}\left( 1+u^{2}\right)
+v^{2}\right] \sqrt{1+u^{2}}}{2\left( 1+2u^{2}\right) v}, &  & \dfrac{%
\partial v}{\partial \varphi }=\sqrt{1-v^{2}}%
\end{array}%
\]%
hence%
\[
\left\vert \frac{\partial (u,v)}{\partial (\varphi ,\psi )}\right\vert =%
\dfrac{\left[ 4u^{2}\left( 1+u^{2}\right) +v^{2}\right] \sqrt{1+u^{2}}\sqrt{%
1-v^{2}}}{2\left( 1+2u^{2}\right) v}
\]%
hence%
\begin{eqnarray*}
\frac{\sin \varphi \cos \psi \sin ^{4}\psi }{\tan ^{2}\psi }d\varphi \,d\psi
&=&\sin \varphi \cos ^{3}\psi \sin ^{2}\psi \,d\varphi \,d\psi  \\
&=&v\left( \frac{4u^{2}\left( 1+u^{2}\right) +v^{2}}{v^{2}}\right) ^{-3/2}%
\left[ 1-\left( \frac{4u^{2}\left( 1+u^{2}\right) +v^{2}}{v^{2}}\right) ^{-1}%
\right] d\varphi \,d\psi  \\
&=&v\,\frac{v^{3}}{\left[ 4u^{2}\left( 1+u^{2}\right) +v^{2}\right] ^{3/2}}\,%
\frac{4u^{2}\left( 1+u^{2}\right) }{4u^{2}\left( 1+u^{2}\right) +v^{2}}%
d\varphi \,d\psi  \\
&=&\frac{v^{4}}{\left[ 4u^{2}\left( 1+u^{2}\right) +v^{2}\right] ^{5/2}}%
\,4u^{2}\left( 1+u^{2}\right) \left\vert \frac{\partial (\varphi ,\psi )}{%
\partial (u,v)}\right\vert du\,dv \\
&=&\frac{v^{4}}{\left[ 4u^{2}\left( 1+u^{2}\right) +v^{2}\right] ^{7/2}}%
\,4u^{2}\sqrt{1+u^{2}}\,\frac{2\left( 1+2u^{2}\right) v}{\sqrt{1-v^{2}}}%
du\,dv.
\end{eqnarray*}%
From%
\[
\dint\limits_{0}^{1}\frac{v^{5}}{\left[ 4u^{2}\left( 1+u^{2}\right) +v^{2}%
\right] ^{7/2}}\,\frac{1}{\sqrt{1-v^{2}}}dv=\frac{4}{15}\frac{1}{\left(
1+2u^{2}\right) ^{6}}\frac{1}{u\sqrt{1+u^{2}}}
\]%
we deduce that the original integral $\zeta _{5}/640$ is equal to%
\begin{eqnarray*}
&&\frac{4}{15\pi ^{2}}\dint\limits_{0}^{\infty }\,\frac{\left[ -2+4\left(
\left( 1+2u^{2}\right) ^{2}-1\right) \right] E\left( i\,u\right) ^{2}}{%
\left( 1+2u^{2}\right) ^{5}}u\,du \\
&&+\frac{8}{15\pi ^{2}}\dint\limits_{0}^{\infty }\,\frac{\left[ 1-2\left(
\left( 1+2u^{2}\right) ^{2}-1\right) +(1+2u^{2})\right] E\left( i\,u\right)
K\left( i\,u\right) }{\left( 1+2u^{2}\right) ^{5}}u\,du \\
&&+\frac{4}{15\pi ^{2}}\dint\limits_{0}^{\infty }\,\frac{\left[ -1+\left(
\left( 1+2u^{2}\right) ^{2}-1\right) -(1+2u^{2})\right] K\left( i\,u\right)
^{2}}{\left( 1+2u^{2}\right) ^{5}}u\,du
\end{eqnarray*}%
which is $(2/5)(\zeta _{4}/256)$. \ \newline

A more complicated but similar argument yields that $\zeta _{6}=\zeta _{4}$.
\ Applying the same procedure for $n=3$, however, is fruitless. The hard
integral in this special case becomes%
\begin{eqnarray*}
&&\dint\limits_{0}^{\pi /2}\,\dint\limits_{0}^{\pi /2}\sin \varphi \sqrt{%
\sin ^{2}\theta \sin ^{2}\varphi +\cos ^{2}\varphi }\,\frac{1}{4\pi }\sin
\varphi \,d\varphi \,d\theta  \\
&=&\dint\limits_{0}^{\pi /2}\,\dint\limits_{0}^{\pi /2}\cos \varphi \sin
^{2}\varphi \sqrt{\sin ^{2}\theta \tan ^{2}\varphi +1}\,\frac{1}{4\pi }%
\,d\theta \,d\varphi  \\
&=&\dint\limits_{0}^{\pi /2}\cos \varphi \sin ^{2}\varphi \,E\left( i\tan
\varphi \right) \frac{1}{4\pi }\,d\varphi  \\
&=&\frac{1}{4\pi }\dint\limits_{0}^{\infty }\,\frac{t^{2}\,E\left(
i\,t\right) }{\left( 1+t^{2}\right) ^{5/2}}\,dt=\frac{1}{96}\zeta _{3}
\end{eqnarray*}%
which bears no obvious resemblance to $(8/3)(\zeta _{4}/256)$. \ Therefore $%
\zeta _{3}=\zeta _{4}$ remains unproven, implying that $\zeta _{4}=3\pi
\,_{3}F_{2}\left( -\tfrac{1}{2},\tfrac{1}{2},\tfrac{3}{2};1,2;1\right) $ is
not yet a theorem. \

\bigskip

\end{document}